\documentclass[12pt,a4paper]{article}
\usepackage{latexsym,amsmath,amsthm,amssymb,mathrsfs}
\usepackage[all]{xy}
\usepackage{textcomp}
\usepackage{graphicx}

\begin{document}

\title{Migration paths saturations in meta-epidemic systems.}

\author{Silvia Motto, Ezio Venturino\\
Dipartimento di Matematica ``Giuseppe Peano'',\\
Universit\`a di Torino, \\via Carlo Alberto 10, 10123 Torino, Italy
}
\date{}

\maketitle

\begin{abstract}
In this paper we consider a simple two-patch model in which a population
affected by a disease can freely move. We assume that the capacity of
the interconnected paths is limited, and thereby influencing the migration
rates. Possible habitat disruptions due to human activities or natural events
are accounted for. The demographic assumptions prevent the ecosystem to
be wiped out, and the disease remains endemic in both populated patches at
a stable equilibrium, but
possibly also with an oscillatory behavior in the case of unidirectional migrations.
Interestingly, if infected cannot migrate, it is possible that one patch becomes disease-free.
This fact could be exploited to keep disease-free at least part of the population.
\end{abstract}

\section{Introduction}

Natural landscapes can become fragmented due to landslides, for instance, or human constructions.
Wild populations can be affected by these events. To understand these phenomena and 
possibly alleviate their negative consequences for the environment, scientists have
developed the concepts of population assembling, \cite{LM}, and
metapopulations, \cite{HG,W97}, which showed that global survival
is possible even if in some patches the populations get extinguished, \cite{W96}.

Diseases represent a common occurrence in nature for individuals and communities.
Ecoepidemiology merges the demographic and epidemiological features of
interacting populations into a single model, see Chapter 7 of \cite{MPV} for an
introduction. In this context also metaecoepidemic models can be considered, \cite{V}.

Epidemics affecting populations living in patchy habitats have been
investigated since quite some time, also in the context of fighting
new emerging diseases, both deterministically, \cite{AvdD,SvdD,W,WR,WZ},
and stochastically, \cite{AP02}.

We consider a very simple one disease-affected population, 2-patch model with migrations.
In \cite{BIV,BCGV} other models of this kind have been introduced. A
specific feature of this contribution lies in the fact that
an upper bound on the migration rates is assumed, as in \cite{ABMV}, but
the restrictive assumption of no vital dynamics used for that similar model
is here removed. Of interest is the assessment of the consequences that
possible paths disturbances have on the whole ecosystem.
The disease is assumed to be recoverable, but both disease transmission and recovery
rates are environment-dependent. The disease tranmission is modeled via
mass action, assuming homogeneous mixing for the population in both patches.
Instead in \cite{ABLN} $n$ patches are assumed, where SIS models with standard incidence
are present in each one. No vital dynamics is however considered.
In \cite{GR} the model is similar to \cite{ABLN}, but contains susceptible recruitments
and a different disease incidence.

The effect of population diffusion on the disease spread is studied in \cite{WM},
investigating what happens if the disease gets eradicated in neighboring patches. Also,
diffusion may or may not help the epidemics to spread, \cite{WR}.

The main reference model is \cite{JW}, where the environment consists of several fragments.
The stability conditions for endemic and disease-free equilibria are established.
In contrast to \cite{JW}, we propose here a saturation
effect on the migrating corridors. The analysis of course shows that the basic equilibria
are the same, i.e.
system disappearance and the endemic equilibrium with both patches populated.
Further, we also try to answer the question of what happens to the ecosystem as a whole
when some disruptions in the interpatch communications occur. This could happen because
migrations are not possible in one direction, or if they require a strenuous effort,
which infected individuals cannot exert.

\section{The Model}\label{chapt:secondo_modello}

A similar model with restricted migrations has been introduced in \cite{ABMV}. But in
contrast to its assumptions stating that the migrations are restricted by the size of the population
of the patch into which the migration occurs, we rather consider here the case in which migrations are restricted
by the size of the available canals. Thus, even if the populations in each patch grow,
only a maximal fixed migration rate can be attained. Mathematically, this is obtained by
using a Holling type II function for modelling the migration rates.

For $k=1,2$, let us denote by $S_k$ the susceptibles and by $I_k$ the infected in each patch.
\begin{eqnarray}\label{model_g}
\dot{S_1}=r_1S_1-\gamma_{1}S_1I_1+\delta_1I_1-m_{21}\frac {S_1}{A+I_1+S_1}+m_{12}\frac {S_2}{A+I_2+S_2},\\ \nonumber
\dot{I_1}=\gamma_{1}S_1I_1-(\delta_1+\mu_1)I_1-n_{21}\frac {I_1}{B+I_1+S_1}+n_{12}\frac {I_2}{B+I_2+S_2},\\ \nonumber
\dot{S_2}=r_2S_2-\gamma_{2}S_2I_2+\delta_2I_2+m_{21}\frac {S_1}{A+I_1+S_1}-m_{12}\frac {S_2}{A+I_2+S_2},\\ \nonumber
\dot{I_2}=\gamma_{2}S_2I_2-(\delta_2+\mu_2)I_2+n_{21}\frac {I_1}{B+I_1+S_1}-n_{12}\frac {I_2}{B+I_2+S_2}.
\end{eqnarray}
The parameters have the following meanings.
By $r_k$ we denote the net reproduction rate of the susceptibles in patch $k$.
Note that we make the strong demographic assumption that only susceptibles give birth,
the disease preventing the infected to reproduce. Further,
$\gamma_k$ denotes the disease contact rate,
$\delta_k$ is the disease recovery rate,
$\mu_k$ is the infected mortality rate in each patch,
$A$ is the half saturation constant for the susceptibles, and $B$ the one for the infected;
finally the migration rates from patch $j$ into patch $i$ are
$m_{ij}$ for the susceptibles and
$n_{ij}$ for the infected.
In fact, e.g. the parameter
$m_{12}$ represents the maximum migration rate possible
for the susceptibles through the canal leading from patch 2 into patch 1.
The last term in the first equation states thus that the higher the population in patch 2,
the smaller the migration rate becomes in view of the saturation of the communication path.
Similarly for the corresponding terms in this and the other equations.

The first equation states that the susceptibles reproduce, and possibly become infected by contagion,
new recruits come into this class also via disease recovery, and then the emigrations and immigrations occur. Similar considerations hold true for the remaining equations.

\begin{figure}
\centering
\includegraphics[scale=0.4]{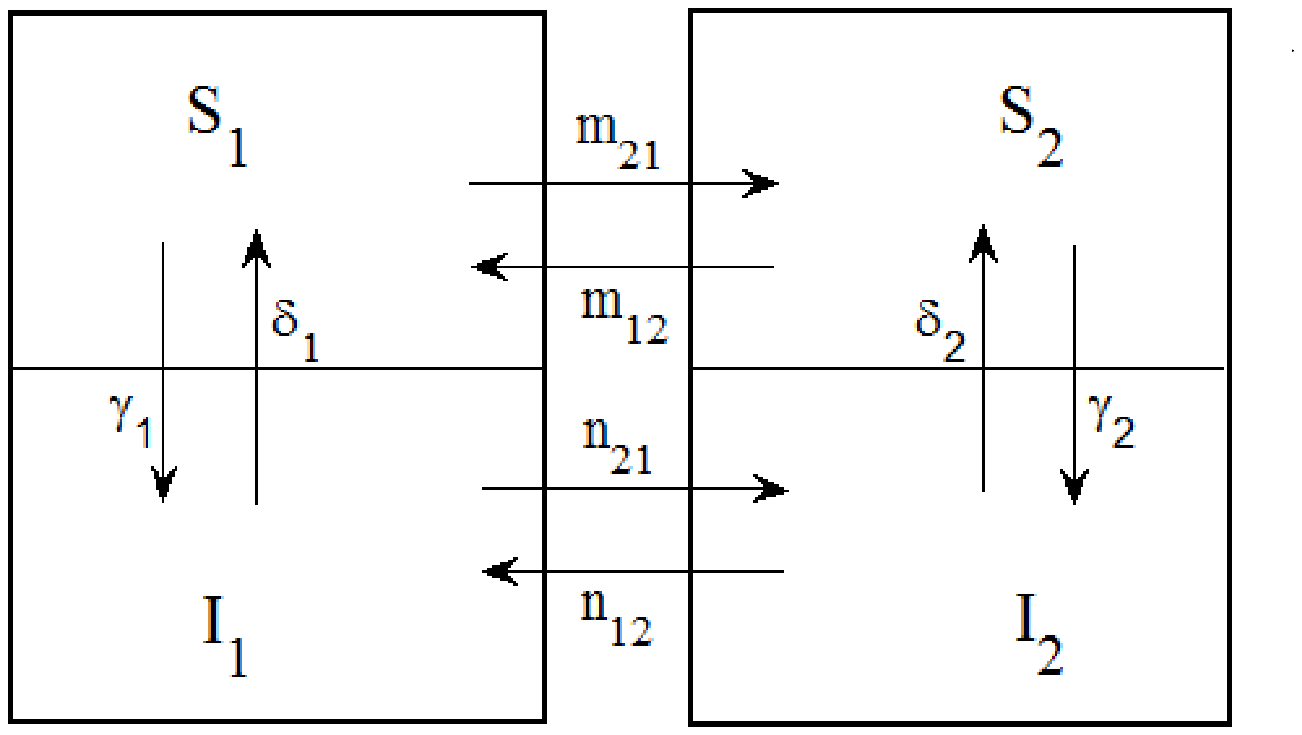}
\label{fig:Illustraz_caso_gen2}
\caption{The system in consideration.}
\end{figure}

\subsection{Equilibria}

The general model admits only two possible equilibria, the trivial state in which the ecosystem vanishes,
and the coexistence state.

To study the latter,
we can eliminate the migration rates by summing the first and third equations of (\ref{model_g}),
as well as the second and fourth one, to obtain respectively
\begin{eqnarray}\nonumber
r_1\tilde{S_1}-\gamma_{1}\tilde{S_1}\tilde{I_1}+\delta_1\tilde{I_1}+r_2\tilde{S_2}-\gamma_{2}\tilde{S_2}\tilde{I_2}+\delta_2\tilde{I_2}=0,\\ \label{eq2piu4}
\gamma_{1}\tilde{S_1}\tilde{I_1}-(\delta_1+\mu_1)\tilde{I_1}+\gamma_{2}\tilde{S_2}\tilde{I_2}-(\delta_2+\mu_2)\tilde{I_2}=0.
\end{eqnarray}
These equations can also be summed, to produce
\begin{equation}
\tilde{S_1}=\frac{-r_2\tilde{S_2}+\mu_1\tilde{I_1}+\mu_2\tilde{I_2}}{r_1}
\end{equation}
which upon substitution into (\ref{eq2piu4}) gives
\begin{equation}
\tilde{S_2}=\frac{r_1((\delta_1+\mu_1)\tilde{I_1}+(\delta_2+\mu_2)\tilde{I_2})-\gamma_1\mu_1{\tilde{I_1}}^2-\gamma_1\mu_1\tilde{I_1}\tilde{I_2}}{r_1\gamma_2\tilde{I_2}-r_2\gamma_1\tilde{I_1}}.
\end{equation}
We need nonnegative populations, therefore some necessary conditions for the
feasibility of the equilibrium with endemic disease and both patches populated follow:
\begin{equation}
\tilde{I_1}>\frac{r_1\gamma_2\tilde{I_2}}{r_2\gamma_1}, \quad
\tilde{S_2} \le \frac{\mu_1\tilde{I_1}+\mu_2\tilde{I_2}}{r_2},  \quad
\tilde{I_2} \ge \frac{\gamma_1\mu_1{\tilde{I_1}}^2-r_1(\delta_1+\mu_1)\tilde{I_1}}{r_1(\delta_2+\mu_2)-\gamma_1\mu_2\tilde{I_1}}
\end{equation}
or the opposite inequalities. The last one, however, leads to an upper bound that must
be explicitly imposed not to be negative. In conclusion, we have the second set of
necessary conditions
\begin{equation}
\tilde{I_1} < \frac{r_1\gamma_2\tilde{I_2}}{r_2\gamma_1}, \quad
\tilde{S_2} \le \frac{\mu_1\tilde{I_1}+\mu_2\tilde{I_2}}{r_2}, \quad
\tilde{I_2} \le \frac{\gamma_1\mu_1{\tilde{I_1}}^2-r_1(\delta_1+\mu_1)\tilde{I_1}}{r_1(\delta_2+\mu_2)-\gamma_1\mu_2\tilde{I_1}},
\end{equation}
supplemented by either one of the two sets of conditions,
$$
\frac {\delta_2 + \mu_2}{\gamma_1 \mu_2} < I_1 \le \frac {\delta_1 + \mu_1}{\gamma_2 \mu_1},
\quad
\frac {\delta_1 + \mu_1}{\gamma_2 \mu_1} \le I_1 < \frac {\delta_2 + \mu_2}{\gamma_1 \mu_2}
$$

\subsection{Stability}

The Jacobian of (\ref{model_g}) is
\begin{equation}\label{Jac}
J=\left[
\begin{array}{cccc}
J_{11}
&-\gamma_1S_1+\delta_1+\hat{\eta_2}S_1&\hat{\theta_1}-\hat{\theta_2}S_2&-\hat{\theta_2}S_2\\
\\
\gamma_1I_1+\hat{\rho_2}I_1& J_{22}
&-\hat{\sigma_2}I_2&\hat{\sigma_1}-\hat{\sigma_2}I_2\\
\\
\hat{\eta_1}-\hat{\eta_2}S_1&-\hat{\eta_2}S_1&J_{33}
&-\gamma_2S_2+\delta_2+\hat{\theta_2}S_2\\
\\
-\hat{\rho_2}I_1&\hat{\rho_1}-\hat{\rho_2}I_1&\gamma_2I_2+\hat{\sigma_2}I_2&J_{44}
\end{array}   \vspace{3pt}
\right]
\end{equation}
with
\begin{eqnarray*}
J_{11}=-\gamma_1I_1-\hat{\eta_1}+\hat{\eta_2}S_1+r_1, \quad
J_{22}=\gamma_1S_1-\delta_1-\hat{\rho_1}+\hat{\rho_2}I_1-\mu_1\\
J_{33}=-\gamma_2I_2-\hat{\theta_1}+\hat{\theta_2}S_2+r_2, \quad
J_{44}=\gamma_2S_2-\delta_2-\hat{\sigma_1}+\hat{\sigma_2}I_2-\mu_2,\\
\hat{\eta_1}=\frac{m_{21}}{A+S_1+I_1}, \quad \hat{\eta_2}=\frac{m_{21}}{{(A+S_1+I_1)}^2},\quad 
\hat{\theta_1}=\frac{m_{12}}{A+S_2+I_2}, \\
\hat{\theta_2}=\frac{m_{12}}{{(A+S_2+I_2)}^2},\quad 
\hat{\rho_1}=\frac{n_{21}}{B+S_1+I_1}, \quad \hat{\rho_2}=\frac{n_{21}}{{(B+S_1+I_1)}^2}\\
\hat{\sigma_1}=\frac{n_{12}}{B+S_2+I_2}, \quad \hat{\sigma_2}=\frac{n_{12}}{{(B+S_2+I_2)}^2}.
\end{eqnarray*}

The origin represents the only case in which the stability study can be performed analytically.
The characteristic equation factorizes, to give $H(\lambda)K(\lambda)=0$, with
$$
K(\lambda)=\lambda^2+(\delta_1+\delta_2+\mu_1+\mu_2+\frac{n_{12}}{B}
+\frac{n_{21}}{B})\lambda+(\delta_2+\mu_2)(\delta_1+\mu_1+
\frac{n_{21}}{B})+(\delta_1+\mu_1)\frac{n_{12}}{B}
$$
and
\begin{equation}\label{H}
H(\lambda)=\lambda^2+(\frac{m_{12}}{A}+\frac{m_{21}}{A}-r_1-r_2)\lambda-\frac{m_{12}r_1}{A}-\frac{m_{21}r_2}{A}+r_1r_2.
\end{equation}
For $K(\lambda)$ all coefficients are positive, so that its roots have both negative real parts.
If we consider the Routh-Hurwitz conditions for $H(\lambda)=0$, we find
$$
\frac{m_{12}}{A}+\frac{m_{21}}{A}-r_1-r_2>0 \quad
-\frac{m_{12}r_1}{A}-\frac{m_{21}r_2}{A}+r_1r_2>0 .
$$
The stability conditions are then
$$
A(r_2+r_1)-m_{21}<m_{12}<Ar_2-\frac{m_{21}r_2}{r_1}.
$$
Eliminating $Ar_2$ from the first and last terms, and observing that $m_{12}>0$ in the last inequality, we get
$$
Ar_1-m_{21}<-\frac{m_{21}r_2}{r_1}  \quad
A-\frac{m_{21}}{r_1}>0 ,
$$
from which $-m_{21}r_2 r_1^{-1}>0$ follows, thus showing that the origin can never be stable.

Through numerical simuations, it can be verified that
indeed the endemic equilibrium can be stably achieved.
This can be accomplished for instance using the
following set of parameter values
\begin{eqnarray*}
r_1=1, \quad r_2=1, \quad \gamma_1=1, \quad \gamma_2=1, \quad \delta_1=0.5, \quad \delta_2=0.5, \quad \mu_1=1,\\
\mu_2=1, \quad m_{12}=1, \quad m_{21}=1, \quad n_{12}=1, \quad n_{21}=1, \quad A=1, \quad B=10.
\end{eqnarray*}

\subsection{Bifurcations}\label{subsec:bif_gen2}
We now show that no Hopf bifurcations can arise at the origin. Since $K(\lambda)$ has roots
with negative real parts, we consider only $H(\lambda)=0$. To have a Hopf bifurcation we need
$$
\frac{m_{12}}{A}+\frac{m_{21}}{A}-r_1-r_2=0, \quad
-\frac{m_{12}r_1}{A}-\frac{m_{21}r_2}{A}+r_1r_2>0 .
$$
Solving for $r_2$ in the first equation, and substituting into the second one,
we have
$$
\Psi(r_1)=
-{r_1}^2+2\frac{m_{21}r_1}{A}-\frac{m_{21}}{A}\left(\frac{m_{12}}{A}+\frac{m_{21}}{A}\right)>0 .
$$
But this condition
can never be satisfied, as the concave parabola $\Psi(r_1)=0$ has the vertex
$\left(m_{21}A^{-1}, -m_{21}m_{12}A^{-2}\right)$, lying in the fourth quadrant.

\section{Unidirectional Migrations}
\label{subsec:no_12_mod2}

\begin{figure}
\centering
\includegraphics[scale=0.4]{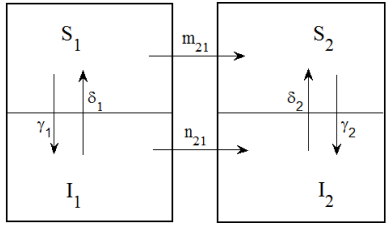}
\caption{The schematic model of unidirectional migrations.}
\label{fig:graf_interr_migr_2_12}
\end{figure}

We now consider the case in which the joining path between the two patches can be
traversed only in one direction. This is by no means restrictive as for instance fish can swim
much more easily downstream in rivers, and sometimes dams and waterfalls prevent them from
returning upstream. The Figure \ref{fig:graf_interr_migr_2_12} describes the situation.

The system (\ref{model_g}) contains now $m_{12}=0$ and $n_{12}=0$.
The system's Jacobian (\ref{Jac}) simplifies accordingly.

\subsection{Equilibria}
\label{subsec:equilibri_no12_mod2}

In this case we have again the origin, and possibly coexistence. But in addition, we find the point
$E_1=(0,0,\tilde{S_2},\tilde{I_2})$ with
\begin{equation}
\tilde{S_2}=\frac{\delta_2+\mu_2}{\gamma_{2}}, \quad
\tilde{I_2}=\frac{r_2(\delta_2+\mu_2)}{\gamma_{2}\mu_2},
\end{equation}
which is always feasible.

We also find the point
$E_2=(\tilde{S_1},0,\tilde{S_2},\tilde{I_2})$ with population values
\begin{equation}
\tilde{S_1}=\frac{m_{21}-r_1A}{r_1}, \quad
\tilde{S_2}=\frac{\delta_2+\mu_2}{\gamma_2}, \quad
\tilde{I_2}=\frac{\gamma_2(m_{21}-r_1A)+r_2(\delta_2+\mu_2)}{\gamma_2\mu_2}.
\end{equation}
It has the following feasibility condition,
\begin{equation}\label{feas_E2}
m_{21}\ge r_1A.
\end{equation}

\subsection{Stability}
At the origin, we find the following eigenvalues,
$$
\lambda_1=r_2 \quad,
\lambda_2=-\delta_2-\mu_2, \quad
\lambda_3=\frac{r_1A-m_{21}}{A} , \quad
\lambda_4=-\frac{\delta_1B+n_{21}+\mu_1B}{B},
$$
from which its unconditional instability is immediate.
Since all the eigenvalues are real, no Hopf bifurcation can arise.

At $E_1$ again it is possible to obtain directly the eigenvalues,
$$
\lambda_{1,2}=\frac{-r_2\delta_2\pm\sqrt{{r_2}^2{\delta_2}^2-4{\mu_2}^2r_2(\mu_2+\delta_2)}}{2\mu_2}, \quad
\lambda_3=\frac{-m_{21}+r_1A}{A}
$$
and $\lambda_4=-[(\delta_1+\mu_1)B+n_{21}]B^{-1}<0$. Since also $\lambda_{1,2}<0$ easily, stability is
regulated by the third eigenvalue, giving
\begin{equation}\label{stab21_caso8_mod2}
r_1A<m_{21}.
\end{equation}
Again, no Hopf bifurcations arise, since
the eigenvalues $\lambda_1$ and $\lambda_2$ can never be purely imaginary, as
the parameters are all positive: $r_2\delta_2 \ne 0$, since $r_2>0$ and $\delta_2>0$.

At $E_2$ once more the eigenvalues are explicitly found,
\begin{eqnarray*}
\lambda_1=\gamma_1\tilde{S_1}-\delta_1-\mu_1-\frac{n_{21}}{B+\tilde{S_1}} , \quad
\lambda_2=-\frac{m_{21}}{A+\tilde{S_1}}+\frac{m_{21}\tilde{S_1}}{{(A+\tilde{S_1})}^2}+r_1 , \\
\lambda_{3,4}=\frac{r_2-\gamma_2\tilde{I_2} \pm \sqrt{(r_2-\gamma_2\tilde{I_2})^2-4\mu_2\gamma_2\tilde{I_2}}}{2}.
\end{eqnarray*}
Using the expression for $\tilde{S_1}$
the second eigenvalue becomes $\lambda_2= r_1(m_{21}-r_1A) m_{21}^{-1}$,
so that its negativity in this case entails
$m_{21}-r_1A<0$,
which contradicts the feasibility condition (\ref{feas_E2}). In
conclusion, $E_2$ is unconditionally unstable.
Also here no Hopf bifurcations arise. Imposing the real part of $\lambda_{3,4}$ to be zero, we find
$r_2-\gamma_2\tilde{I_2}=0$ which explicitly becomes
\begin{displaymath}
-\frac{\gamma_2(m_{21}-r_1A)+r_2\delta_2}{\mu_2}=0,
\end{displaymath}
which cannot be satisfied in view of the feasibility condition (\ref{feas_E2}).

In this model we can study the coexistence equilibrium, because the Jacobian becomes a lower triangular
matrix. The characteristic equation then factorizes accordingly, to give the quadratic equation
\begin{equation}
{\lambda}^2+a_1\lambda+a_0=0
\end{equation}
with
\begin{displaymath}
\begin{array}{l}
a_1=\gamma_1\tilde{I_1}+\tilde{\eta_1}-\tilde{\eta_2}\tilde{S_1}-r_1-\gamma_1\tilde{S_1}+\delta_1+\mu_1+\tilde{\rho_1}-\tilde{\rho_2}\tilde{I_1} \\
a_0=(-\gamma_1\tilde{I_1}-\tilde{\eta_1}+r_1)(-\delta_1-\mu_1-\tilde{\rho_1}+\tilde{\rho_2}\tilde{I_1}) 
+\gamma_1\tilde{S_1}\tilde{I_1}(\tilde{\rho_2}-\tilde{\eta_2})\\
+\tilde{\eta_2}\tilde{S_1}(\gamma_1\tilde{S_1}-\delta_1-\mu_1
-\tilde{\rho_1})-\delta_1\tilde{I_1}(\gamma_1+\tilde{\rho_2}).
\end{array}
\end{displaymath}
To have roots with negative real parts we then need both these coefficients positive,
\begin{displaymath}
a_1>0, \quad 
a_0>0.
\end{displaymath}
The remaining eigenvalues are evaluated explicitly,
\begin{displaymath}
\lambda_{3,4}=\frac{k \pm \sqrt{k^2+4(r_2\delta_2+r_2\mu_2-r_2\gamma_2\tilde{S_2}-\mu_2\gamma_2\tilde{I_2})}}{2},
\end{displaymath}
where $k=-\gamma_2\tilde{I_2}-\delta_2-\mu_2+\gamma_2\tilde{S_2}+r_2$. They are both real
and negative if
$k<0$ and $r_2\delta_2+r_2\mu_2-r_2\gamma_2\tilde{S_2}-\mu_2\gamma_2\tilde{I_2}<0$.
In summary, the stability conditions for the coexistence equilibrium are
\begin{equation}\label{stab16}
a_1>0, \quad
a_0>0, \quad
\gamma_2\tilde{S_2}+r_2<\gamma_2\tilde{I_2}+\delta_2+\mu_2, \quad
r_2\delta_2+r_2\mu_2<r_2\gamma_2\tilde{S_2}+\mu_2\gamma_2\tilde{I_2}.
\end{equation}

In principle Hopf bifurcations could arise in this situation, whenever either one of the two
sets of conditions holds,
\begin{equation}\label{condbif1}
a_1=0, \quad a_0>0, \quad k<0, \quad h<0,
\end{equation}
or
\begin{equation}\label{condbif2}
a_1>0, \quad a_0>0, \quad k=0, \quad h<0.
\end{equation}

\section{Infected do not Migrate}
\label{sec:terzo_modello}

In this case we assume that migrations entail an effort, which
is too strenuous for infected to exert, so that they are prevented
from changing the patch in which they live.
We need to set $n_{21}=n_{12}=0$ and dropping also the populations
$I_1$ and $I_2$ from the migration terms.

\begin{figure}
\centering
\includegraphics[scale=0.4]{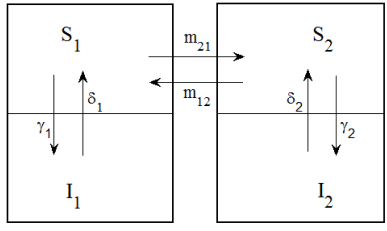}
\caption{The system where the infected are prevented from migrating.}
\label{fig:infetti_non_migrano}
\end{figure}

Pictorially, the system is illustrated in Figure \ref{fig:infetti_non_migrano}.

\subsection{Equilibria}
\label{subsec:equilibri_gen3}

In addition to the origin, we find also the following two pairs of equilibria,
$Z_1^{\pm}=(\tilde{S_1},\tilde{I_1},\tilde{S_2}^{\pm},0)$,
$Z_2^{\pm}=(\tilde{S_1}^{\pm},0,\tilde{S_2},\tilde{I_2})$.
At $Z_1$ the population values can be explicitly calculated,
$$
\tilde{S_1}=\frac{\delta_1+\mu_1}{\gamma_1},\quad
\tilde{I_1}=\tilde{S_1}\frac{r_1}{\mu_1}+\frac{r_2\tilde{S_2}}{\mu_1},\quad
\tilde{S_2}^{\pm}=\frac{\ell \pm \sqrt{\ell^2-4r_2m_{21}A\tilde{S_1}(A+\tilde{S_1})}}{2r_2(A+\tilde{S_1})}
$$
where $\ell=m_{12}A+m_{12}\tilde{S_1}-r_2A^2-r_2A\tilde{S_1}-m_{21}\tilde{S_1}$.
These points are feasible if and only if
\begin{equation}\label{amm_12_3bis}
m_{21}\tilde{S_1} +r_2A^2+r_2A\tilde{S_1} <  m_{12}A+m_{12}\tilde{S_1}.
\end{equation}
The inequality is strict since for $\ell =0$ the corresponding quadratic equation has purely
imaginary solutions.
We also find at $Z_2$ the populations
\begin{eqnarray*}
\tilde{S_2}=\frac{\delta_2+\mu_2}{\gamma_2}, \quad
\tilde{I_2}=r_2\frac{\delta_2+\mu_2}{\gamma_2\mu_2}+\frac{r_1\tilde{S_1}}{\mu_2},\quad
\tilde S_1^{\pm}=\frac{h \pm \sqrt{h^2-4r_1m_{12}A\tilde{S_2}(A+\tilde{S_2})}}{2r_1(A+\tilde{S_2})}
\end{eqnarray*}
with $h=m_{21}A+m_{21}\tilde{S_2}-m_{12}\tilde{S_2}-r_1A^2-r_1A\tilde{S_2}$.
Once again, noting the strict inequality, these equilibria are feasible if and only if
\begin{equation}   \label{amm_14_3bis}
m_{12}\tilde{S_2}+r_1A^2+r_1A\tilde{S_2} < m_{21}A+m_{21}\tilde{S_2} .
\end{equation}

Also the endemic coexistence equilibrium can be analytically evaluated,
\begin{eqnarray*}
\tilde{S_1}=\frac{\delta_1+\mu_1}{\gamma_{1}}, \quad
\tilde{I_2}=\frac{r_1\gamma_2(\delta_1+\mu_1)+r_2\gamma_1(\delta_2+\mu_2)
-\tilde{I_1}\gamma_1\gamma_2\mu_1
}{\gamma_1\gamma_2\mu_2}, \\
\tilde{S_2}=\frac{\delta_2+\mu_2}{\gamma_{2}},\quad
\tilde{I_1}
= \frac 1{\gamma_{1}\tilde{S_1}-\delta_1}
\left[r_1\tilde{S_1}-m_{21}\frac{\tilde{S_1}}{A+\tilde{S_1}}+m_{12}\frac{\tilde S_{2}}{A+\tilde{S_2}}\right]
\end{eqnarray*}

For feasibility, note that the denominator in the expression for $\tilde I_1$ reduces
to $\mu_1$. Then $I_1 \ge 0$ gives
\begin{equation}   \label{feas_coex1}
r_1\tilde{S_1}+ m_{12}\frac{\tilde S_{2}}{A+\tilde{S_2}} \ge 
m_{21}\frac{\tilde{S_1}}{A+\tilde{S_1}}.
\end{equation}
We need also $I_2 \ge 0$ i.e.
\begin{equation}\label{amm_coex}
\tilde{I_1}\le \frac{r_1\gamma_2(\delta_1+\mu_1)+r_2\gamma_1(\delta_2+\mu_2)}{\gamma_1\gamma_2\mu_1}
\end{equation}
which can be recast in the following form
\begin{equation}   \label{feas_coex2}
r_2\tilde{S_2}+m_{21}\frac{\tilde{S_1}}{A+\tilde{S_1}} \ge m_{12}\frac{\tilde S_{2}}{A+\tilde{S_2}}.
\end{equation}

\subsection{Stability}
At the origin, the Jacobian has the explicit eigenvalues
$\lambda_1=-\delta_1-\mu_1$, $\lambda_2=-\delta_2-\mu_2$
and the roots of the quadratic (\ref{H}) so that the same analysis
carries out also in this case, showing the inconditionate
instability of this equilibrium point.

At $Z_1$ one eigenvalue is $\lambda_1=\gamma_2\tilde{S_2}-\delta_2-\mu_2$.
The remaining ones are the roots of the cubic equation
$\lambda^3+\hat{p_2}\lambda^2+\hat{p_1}\lambda+\hat{p_0}=0$
with
\begin{eqnarray*}
\hat{p_2}=\gamma_1\tilde{I_1}-r_2-r_1+\tilde{\eta_1}-\tilde{\eta_2}\tilde{S_1}+\tilde{\theta_1}-\tilde{\theta_2}\tilde{S_2}-\gamma_1\tilde{S_1}+\delta_1+\mu_1, \\
\hat{p_1}=\gamma_1\tilde{I_1}(\gamma_1\tilde{S_1}-\delta_1)+(\tilde{\theta_1}-\tilde{\theta_2}\tilde{S_2})(\gamma_1\tilde{I_1}-r_1-\gamma_1\tilde{S_1}+\delta_1 +\mu_1) \\
          +(\tilde{\eta_1}\tilde{S_1}+r_1-\gamma_1\tilde{I_1}-\tilde{\eta_1})(r_2+\gamma_1\tilde{S_1}-\delta_1-\mu_1)
          +r_2(\gamma_1\tilde{S_1}-\delta_1-\mu_1), \\
\hat{p_0}=\gamma_1\tilde{I_1}(\gamma_1\tilde{S_1}-\delta_1)(\tilde{\theta_1}-\tilde{\theta_2}\tilde{S_2}-r_2)-(\gamma_1\tilde{S_1}-\delta_1-\mu_1)[r_2(\tilde{\eta_2}\tilde{S_1} \\
          +r_1-\gamma_1\tilde{I_1}-\tilde{\eta_1})+(\gamma_1\tilde{I_1}-r_1)(\tilde{\theta_1}-\tilde{\theta_2}\tilde{S_2})].
\end{eqnarray*}

The Routh-Hurwitz conditions combined with negativity of the explicit eigenvalue ensure then stability:
\begin{equation}\label{stab_gen3_caso12}
\gamma_2\tilde{S_2}<\delta_2+\mu_2, \quad
\hat{p_0}>0, \quad
\hat{p_2}>0, \quad
\hat{p_2}\hat{p_1}>\hat{p_0}.
\end{equation}

A similar situation arises for $Z_2$, one eigenvalue is found analytically,
$\lambda_1=\gamma_1\tilde{S_1}-\delta_1-\mu_1$
and the cubic equation $\lambda^3+\hat{q_2}\lambda^2+\hat{q_1}\lambda+\hat{q_0}=0$ with coefficients
\begin{eqnarray*}
\hat{q_2}=-r_2-r_1+\gamma_2\tilde{I_2}+\tilde{\theta_1}-\tilde{\theta_2}\tilde{S_2}+\tilde{\eta_1}-\tilde{\eta_2}\tilde{S_1}-\gamma_2\tilde{S_2}+\delta_2+\mu_2 \\
\hat{q_1}=-\gamma_2\tilde{I_2}(-\gamma_2\tilde{S_2}+\delta_2)+(-\tilde{\eta_1}+\tilde{\eta_2}\tilde{S_1})(-\gamma_2\tilde{I_2}+r_2+\gamma_2\tilde{S_2}-\delta_2 -\mu_2) \\
          +(-\gamma_2\tilde{I_2}-\tilde{\theta_1}+\tilde{\theta_1}\tilde{S_2}+r_2)(r_1+\gamma_2\tilde{S_2}-\delta_2-\mu_2)+ 
          +r_1(\gamma_2\tilde{S_2}-\delta_2-\mu_2) \\
\hat{q_0}=\gamma_2\tilde{I_2}(-\gamma_2\tilde{S_2}+\delta_2)(-\tilde{\eta_1}+\tilde{\eta_2}\tilde{S_1}+r_1)-(\gamma_2\tilde{S_2}-\delta_2-\mu_2)[r_1(-\gamma_2\tilde{I_2}+ \\
          -\tilde{\theta_1}+\tilde{\theta_2}\tilde{S_2}+r_2)+(-\gamma_2\tilde{I_2}+r_2)(-\tilde{\eta_1}+\tilde{\eta_2}\tilde{S_1})]
\end{eqnarray*}
for which the stability criterion becomes
\begin{equation}\label{stab_gen3_caso14}
\gamma_1\tilde{S_1}<\delta_1+\mu_1, \quad
\hat{q_0}>0 , \quad
\hat{q_2}>0 , \quad
\hat{q_2}\hat{q_1}>\hat{q_0} 
\end{equation}

\paragraph{Remark 1.}

No bifurcations can arise here near the origin. The proof of this statement is exactly the same as the one
carried out in Section \ref{subsec:bif_gen2}.

The stability conditions for both equilibria $Z_1$ and $Z_2$ are nonempty, as can be easily shown numerically
using respectively the following sets of parameters
\begin{eqnarray*}
r_1=1, \quad r_2=0.8, \quad \gamma_1=0.5, \quad \gamma_2=1, \quad \delta_1=1, \quad \delta_2=4,\\
\mu_1=1, \quad \mu_2=2, \quad m_{21}=2, \quad m_{12}=10, \quad A=5.
\end{eqnarray*}
and
\begin{eqnarray*}
r_1=0.8, \quad r_2=1, \quad \gamma_1=1, \quad \gamma_2=0.5, \quad \delta_1=4, \quad \delta_2=1,\\
\mu_1=2, \quad \mu_2=1, \quad m_{21}=9, \quad m_{12}=2, \quad A=5.
\end{eqnarray*}

The endemic equilibrium with all patches populated can numerically be shown to be attained for instance
for the parameter values
\begin{eqnarray*}
r_1=1, \quad r_2=1, \quad \gamma_1=1, \quad \gamma_2=1, \quad \delta_1=0.5, \quad \delta_2=0.5,\\
\mu_1=1, \quad \mu_2=1, \quad m_{12}=1, \quad m_{21}=1, \quad A=1.
\end{eqnarray*}

\section{Biological Interpretation}

For the general model, the system can never be wiped out, since the origin is unconditionally unstable.
The ecosystem thrives with a nonvanishing population and an endemic state of the disease in both patches at stable levels, for certain
parameter ranges.

These results hold true also for the particular case in which migrations back into patch 1 are forbidden. But in such case new possible
equilibria arise, in which patch 1 is depleted, or in which only the susceptible population survives. But the latter equilibrium is
never stable.
The equilibrium with patch 1 empty is stable if the reproductive rate in that patch is low enough, or better, if the emigration rate
is sufficiently high, compare (\ref{stab21_caso8_mod2}).
For the equilibrium with endemic disease and both patches populated, stability conditions have been derived, and the presence
of a regime of possible oscillatory behavior has been highlighted.

If the infected do not migrate, once again the ecosystem is guaranteed to survive, as the origin is always unstable.
The equilibria $Z_1$ and $Z_2$ are interesting, as in them one patch becomes disease-free. This is a result that potentially could be
exploited by the manager of wild parks, to preserve at least part of a population from an epidemics.
In order to control the disease at least in part of the environment
it therefore appears to be better to preserve population movements in both directions,
by preventing the infected to migrate, than
to impose unidirectional migrations for both classes of individuals.

\subsubsection*{Acknowledgments.} The authors have been partially supported by the project ``Metodi numerici
in teoria delle popolazioni'' of the Dipartimento di Matematica ``Giuseppe Peano''.


\begin{thebibliography}{4}
\bibitem{ABMV} Aimar, V., Borlengo, S., Motto, S., Venturino, E.,
A meta-epidemic model with steady state demographics and migrations saturation,
AIP Conf. Proc. 1479, ICNAAM 2012, T. Simos, G. Psihoylos, Ch. Tsitouras, Z. Anastassi (Editors),
1311--1314 (2012); doi: 10.1063/1.4756396

\bibitem{ABLN} Allen, L. J. S., Bolker, B. M., Lou, Y., Nevai, A. L.,
Asymptotic Profiles of the Steady States
for an SIS Epidemic Patch Model, SIAM J. Appl. Math. 67, 1283--1309, (2007)

\bibitem{AvdD} Arino, J., van den Driessche , P.,
Disease spread in metapopulations, Nonlinear
Dynamics and Evolution Equations, Fields Inst. Commun. 48, H. Brunner, X. O. Zhao, and
X. Zou, eds., AMS, Providence, RI, 1--13, (2006).

\bibitem{AP02} Arrigoni, F., Pugliese, A.,
Limits of a multi-patch SIS epidemic model, J. Math. Biol. 45, 419--440 (2002)

\bibitem{BIV} Barengo, M., Iennaco, I., Venturino, E.,
A simple meta-epidemic model,
Proceedings of the 12th International Conference on Computational and 
Mathematical Methods in Science and Engineering, CMMSE 2012, J. Vigo-Aguiar, A.P. Buslaev,
A. Cordero, M. Demiralp, I.P. Hamilton,
E. Jeannot, V.V. Kozlov, M.T. Monteiro, J.J. Moreno, J.C. Reboredo, P. Schwerdtfeger,
N. Stollenwerk, J.R. Torregrosa, E. Venturino, J. Whiteman (Editors)
La Manga, Spain, July 2nd-5th, 2012, v. 1, p. 122--133 (2012)

\bibitem {BCGV} Bianco, F., Cagliero, E., Gastelurrutia, M., Venturino, E.,
Metaecoepidemics with migration of and disease
in the predators, Proceedings of the 11th International Conference on Computational and 
Mathematical Methods in Science and Engineering, CMMSE 2011,
J. Vigo Aguiar, R. Cortina, S. Gray, J.M. Ferradiz, A. Fernandez, I. Hamilton, J.A. Lopez-Ramos,
F. de Oliveira, R. Steinwandt, E. Venturino, J. Whiteman, B. Wade (Editors),
Benidorm, Spain, June 26th-30th, v. 1, 204--223 (2011)

\bibitem{GR} Gao, D., Ruan, S.,
An SIS patch model with variable transmission coefficients,
Mathematical Biosciences 232, 110--115 (2011)

\bibitem {HG} Hanski, I., Gilpin, M.,
Metapopulation biology: ecology, genetics and evolution,
London: Academic Press (1997)

\bibitem{JW} Jin, Y., Wang, W.,
The effect of population dispersal on the spread of a disease, J. Math.
Anal. Appl. 308, 343--364 (2005)

\bibitem{LM} Lloyd, A., May, R. M.,
Spatial heterogeneity in epidemic models, J. Theoret. Biol. 179,
1--11 (1996)

\bibitem{MPV} Malchow, H., Petrovskii, S., Venturino, E.,
Spatiotemporal patterns in Ecology and Epidemiology.
CRC, Boca Raton (2008)

\bibitem{SvdD} Salmani, M., van den Driessche, P.,
A model for disease transmission in a patchy
environment, Discrete Contin. Dynam. Systems Ser. B 6, 185--202 (2006)

\bibitem{V} Venturino, E.,
Simple metaecoepidemic models, Bull. Math. Biol. 73, 917--950 (2011)

\bibitem{W} Wang, W.,
Population dispersal and disease spread,
Discrete Contin. Dynam. Systems Ser. B 4, 797--804 (2004)

\bibitem{WM} Wang, W., Mulone, G.,
Threshold of disease transmission in a patch environment, J.
Math. Anal. Appl. 285, 321--335 (2003)

\bibitem{WR} Wang, W.. Ruan, S.,
Simulating the SARS outbreak in Beijing with limited data, 
J. Theoret. Biol. 227, 369--379  (2004)

\bibitem{WZ} Wang, W., Zhao, X.-Q.,
An epidemic model in a patchy environment, Math. Biosci. 190,
97--112 (2004)

\bibitem {W96} Wiens, J. A.,
Wildlife in patchy environments: metapopulations, mosaics, and management,
in D. R. McCullough (Ed.) Metapopulations and Wildlife Conservation.
Island Press, Washington 53--84 (1996)

\bibitem {W97} Wiens, J. A.,
Metapopulation dynamics and landscape ecology, in I. A. Hanski, M. E. Gilpin
(Ed.s), 43--62. Academic Press, San Diego, (1997)

\end{thebibliography}
\end{document}